\documentclass[12pt]{article}
\usepackage{amsfonts}

\usepackage{latexsym}
\usepackage{amssymb}
\usepackage{epsfig}
\usepackage{amsmath}
\usepackage{amsthm}
\usepackage[mathscr]{eucal}
\usepackage{subfigure}
\usepackage{graphicx}

\newtheorem{Lemma}{Lemma}

\newtheorem{Theorem}[Lemma]{Theorem}

\newtheorem{Definition}{Definition}

\renewcommand{\qed}{\hfill{\ \ \rule{2mm}{2mm}} \vspace{0.2in}}

\begin{document}

\title{Duality between star and plus connected components in percolation}
\author{ \textbf{Ghurumuruhan Ganesan}
\thanks{E-Mail: \texttt{gganesan82@gmail.com} } \\
%EndAName
\ \\
NISER, Bhubaneshwar, India}
\date{}
\maketitle

\begin{abstract}

%In this paper, we study the structure of left-right crossings
%of the random geometric graph \(G = G(n,r_n)\) of \(n\) nodes
%uniformly distributed in \(S = [0,1]^2\) with \(r_n = \epsilon\sqrt{\frac{\log{n}}{n}}\)
%for some \(\epsilon > 0.\) Tiling \(S\) horizontally and
%vertically into rectangles of length \(1\) and width \(Mr_n,\) we
%show that each rectangle has a left-right crossing of edges with
%high probability if \(M\) is sufficiently large.
%We call the resulting subgraph to be a ``backbone" of \(G.\)

%The techniques we use to construct the backbone has quite a few applications.
%As a first, we show that the diameter of second largest component in \(G\)
%is \(O(1)\) with high probability. Secondly,

Tile \(\mathbb{R}^2\) into disjoint unit squares \(\{S_k\}_{k \geq 0}\) with the origin being the centre of \(S_0\)
and say that \(S_i\) and \(S_j\) are star adjacent if they share a corner and plus adjacent if they share an edge.
Every square is either vacant or occupied. Star and plus connected components containing the origin have been previously
studied using unicoherence and interface graphs. In this paper, we use the structure of the outermost boundaries derived in
Ganesan (2015) to alternately obtain duality between star and plus connected components in the following sense:
There is a plus connected cycle of vacant squares attached to surrounding the finite star connected component
containing the origin. There is a star connected cycle of vacant squares attached to and surrounding
the finite plus connected component containing the origin.

\vspace{0.1in} \noindent \textbf{Key words:} Star and plus connected components, duality.

\vspace{0.1in} \noindent \textbf{AMS 2000 Subject Classification:} Primary:
60J10, 60K35; Secondary: 60C05, 62E10, 90B15, 91D30.
\end{abstract}

\bigskip

\renewcommand{\theequation}{\thesection.\arabic{equation}}
\setcounter{equation}{0}
\section{Introduction} \label{intro}

Tile \(\mathbb{R}^2\) into disjoint unit squares \(\{S_k\}_{k \geq 0}\) with origin being the centre of~\(S_0.\) We say \(S_1\) and \(S_2\) are \emph{adjacent} or \emph{star-adjacent} if they share a corner between them. We say that squares \(S_1\) and \(S_2\) are \emph{plus-adjacent}, if they share an edge between them. Here we follow the notation of Penrose (2003). Every square is assigned one of the two states: occupied or vacant.

Let \(C(0)\) denote the star-connected occupied component containing the origin and throughout we assume that \(C(0)\) is finite. Thus if \(S_0\) is vacant then \(C(0) = \emptyset.\) Else \(S_0 \in C(0)\) and if \(S_1, S_2 \in C(0)\) there exists a sequence of distinct occupied squares \((Y_1,Y_2,...,Y_t)\) all belonging to \(C(0),\) such that \(Y_i\) is adjacent to \(Y_{i+1}\) for all \(i\) and \(Y_1 = S_1\) and \(Y_t = S_2.\) Let \(G_C\) be the graph with vertex set being the set of all corners of the squares \(\{S_k\}_k\) in \(C(0)\) and edge set consisting of the edges of the squares \(\{S_k\}_k\) in~\(C(0).\)

Two vertices \(u\) and \(v\) are said to be adjacent in \(G_C\) if they share an edge between them. We say that an edge \(e\) in \(G_C\) is adjacent to square \(S_k\) if it is one of the edges of \(S_k.\) We say that \(e\) is a \emph{boundary edge} if it is adjacent to a vacant square and is also adjacent to an occupied square. A \emph{path} \(P\) in \(G_C\) is a sequence of distinct vertices \((u_0,u_1,...,u_t)\) such that \(u_i\) and \(u_{i+1}\) are adjacent for every \(i.\)  A \emph{cycle} \(C\) in \(G_C\) is a sequence of distinct vertices \((v_0,v_1,...,v_m,v_0)\) starting and ending at the same point such that \(v_i\) is adjacent to \(v_{i+1}\) for all \(0 \leq i \leq m-1\) and \(v_m\) is adjacent to \(v_0.\) A \emph{circuit} \(C'\) in \(G_C\) is a sequence of vertices \((w_0,w_1,...,w_r,w_0)\) starting and ending at the same point such that \(w_i\) is adjacent to \(w_{i+1}\) for all \(0 \leq i \leq r-1,\) \(w_r\) is adjacent to \(w_0\) and no edge is repeated in \(C'.\) Thus vertices may be repeated in circuits and for more related definitions, we refer to Chapter~1, Bollobas (2001).

Any cycle \(C\) divides the plane \(\mathbb{R}^2\) into two disjoint connected regions. As in Bollobas and Riordan (2006), we denote the bounded region to be the \emph{interior} of \(C\) and the unbounded region to be the \emph{exterior} of \(C.\) We have the following definition.
\begin{Definition} We say that edge \(e\) in \(G_C\) is an \emph{outermost boundary} edge of the component \(C(0)\) if the following holds true for every cycle \(C\) in \(G_C:\) either \(e\) is an edge in \(C\) or \(e\) belongs to the exterior of \(C.\)

We define the outermost boundary \(\partial _0\) of \(C(0)\) to be the set of all outermost boundary edges of \(G_C.\)
\end{Definition}
Thus outermost boundary edges cannot be contained in the interior of any cycle in \(G_C.\) In Theorem~1 of Ganesan (2015), we obtained that the outermost boundary for the finite star connected component containing the origin is a connected union of cycles with mutually disjoint interiors. Similary in Theorem~2 of Ganesan (2015), we obtained that the outermost boundary for the finite plus connected component is a single cycle. Using these results, we obtain duality between star and plus connected components.

\subsection*{Duality}
To study duality between star and plus connected components, we first define \(S-\)cycles. Let \(\{J_i\}_{1 \leq i \leq m}\) be a set of squares in \(\cup\{S_k\}_k.\) We say that the sequence \(L = (J_1,...,J_m)\) is a plus connected \(S-\)path if for each \(1 \leq i \leq m-1,\) we have \(J_i\) is plus adjacent to \(J_{i+1}.\) We say that \(L^+\) is a plus connected \(S-\)cycle if in addition \(J_m\) is also  plus adjacent to \(J_1.\) Let \(\Lambda_0\) denote the set of all vacant squares that share a corner with some occupied square in the star connected component \(C(0).\) We have the following result.
\begin{Theorem}\label{thm4} Suppose \(C(0)\) is finite. There exists a unique plus connected \(S-\)cycle \(G_{out} = (Y_1,...,Y_t)\) with the following properties:\\
(i) For every \(i, 1 \leq i \leq t,\) the square \(Y_i \in \cup_k \{S_k\}\) is in \(\Lambda_0.\)\\
(ii) The outermost boundary of \(G_{out}\) is a single cycle \(\partial_G.\)\\
(iii) If \(S_k \in C(0)\cup \Lambda_0,\) then either \(S_k\) is in \(G_{out}\) or \(S_k\) is contained in the interior of \(\partial_G.\)
\end{Theorem}
Thus the corresponding sequence of squares \((Y_{1},...,Y_{t})\) form a plus connected cycle of vacant squares containing all squares of \(C(0)\) in the interior.

To study plus connected components, we define star connected \(S-\)cycles analogously as in plus connected case. As before, let \(\{H_i\}_{1 \leq i \leq n}\) be a set of squares in \(\cup\{S_k\}_k.\) We say that the sequence \(L' = (H_1,...,H_n)\) is a star connected \(S-\)path is for each \(1 \leq i \leq n-1,\) we have \(H_i\) is star adjacent to \(H_{i+1}.\) We say that \(L'\) is a star connected \(S-\)cycle if in addition \(H_n\) is also star adjacent to \(H_1.\) Let \(\Lambda^+_0\) denote the set of all vacant squares that share a corner with some occupied square in the plus connected component \(C^+(0).\) We have the following result.
\begin{Theorem}\label{thm5} Suppose \(C^+(0)\) is finite. There exists a unique star connected \(S-\)cycle \(H_{out} = (U_1,...,U_q)\) with the following properties:\\
(i) For every \(i, 1 \leq i \leq q,\) the square \(U_i \in \cup_k \{S_k\}\) is in \(\Lambda^+_0.\)\\
(ii) The outermost boundary of \(H_{out}\) is a single cycle \(\partial_H.\)\\ %where each \(f_i\) is an edge attached to a vacant square in \(\Lambda_0^+.\)\\
(iii) If \(S_k \in C^+(0)\cup \Lambda^+_0,\) then either \(S_k\) is in \(H_{out}\) or \(S_k\) is contained in the interior of  \(\partial_H.\)
\end{Theorem}
Thus the corresponding sequence of squares \((U_{1},...,U_{q})\) form a star connected cycle of vacant squares containing all squares of \(C(0)\) in the interior. We remark that the proof techniques of the above two results can also be used to study star and plus connected left right and top down crossings in rectangles. For more material we refer to Penrose (2003) and Bollobas and Riordan (2006).

The paper is organized as follows: We prove Theorems~\ref{thm4} and~\ref{thm5} in Sections~\ref{pf4} and~\ref{pf5}, respectively.

%be the square in \(C(0)\) whose centre has the largest \(x\) and

%To see for circuit...

%now to see that every edge is outer most boundary....seems ok...and to write crflly...fr our bnft etc... also to obtain circuit...and write carefully that cycle graph is a tree etc...for our bnft etc..and to write uniqueness by construction the outermost boundary is unique?? by construction shud be unique...suppose there exists a different set of cycles ?? adn argue briefly?? to see dtls crflly..for our bnt etc........

%\begin{figure}[tbp]
%\centering
%\includegraphics[width=3.5in, trim= 0 250 0 275, clip=true]{ckt_mg.eps}
%\caption{Merging cycle \(ABCDA\) with the segment \(AEC\)}
%\label{merg_cyc}
%\end{figure}

\section{Proof of Theorem~\ref{thm4}}\label{pf4}

%(AND TO SEE ALL DETAILS ETC AGAIN CRFLLY..ETC..IMP...

\emph{Proof of Theorem~\ref{thm4}}: To see that such a cycle \(G_{out}\) exists, we let \(\partial_0\) denote the outermost boundary for \(C(0).\) From Theorem~1 of Ganesan (2015), we know that \(\partial_0 = \cup_{1 \leq i \leq n}C_i\) is a connected union of cycles \(\{C_i\}_i\) with mutually disjoint interiors and with the property that \(C_i\) and \(C_j\) intersect at at most one point for disinct \(i\) and \(j.\) Place a unit square on each vertex of \(\partial_0\) and call the union of squares as \(C_V(\partial_0).\) Here the unit square has edges parallel to the \(X-\) and \(Y-\)axes.

By construction, the union \(C_V(\partial_0)\) is a plus connected component. Therefore using Theorem~2 of Ganesan (2015), we have that the outermost boundary \(\partial_V(\partial_0)\) of \(C_V(\partial_0)\) is a single cycle containing all squares of \(C_V(\partial_0)\) in its interior. In particular, \(\partial_0\) is contained in its interior. Here, the outermost boundary \(\partial_V(\partial_0)\) is obtained as follows. Each square in \(C_V(\partial_0)\) is labelled \(1\) and each square sharing a vertex with a square in \(C_V(\partial_0)\) and not belonging to \(C_V(\partial_0)\) is labelled \(0.\) We then apply Theorem~2 of Ganesan (2015) with label \(1\) squares as being occupied and label \(0\) squares as being vacant.

Let \(G_B\) be the graph with vertex set as the \emph{centres} of the squares \(\{S_k\}_k\) and edges drawn between centres of plus adjacent squares in \(\{S_k\}_k.\) The cycle \(\partial_V(\partial_0) = (z_1,...,z_t)\) in \(G_{B}\) satisfies the following properties:\\
(a) For every \(i, 1 \leq i \leq t,\) the square \(Y_i \in \cup_k \{S_k\}\) with centre \(z_i\) is in \(\Lambda_0.\)\\
(b) If \(S_k \in C(0)\cup \Lambda_0,\) then either the centre of \(S_k\) is in \(\partial_V(\partial_0)\) or \(S_k\) is contained in the interior of \(\partial_V(\partial_0).\)

%(c) For every \(2 \leq i \leq t-1,\) the square \(Y_i\) is plus adjacent to only \(Y_{i-1}\) and \(Y_{i+1}\) and \(Y_t = Y_1\) is only plus adjacent to \(Y_{t-1}\) and \(Y_2.\)

Assuming (a)-(b) for the moment, we have that the sequence of squares \(L = (Y_1,...,Y_t)\) is the desired cycle~\(G_{out}\) in the statement of the theorem. Indeed, by construction it satisfies (i). Also \(L\) is a plus connected component and by Theorem~2 of Ganesan (2015), the outermost boundary \(\partial_L(L)\) of \(L\) is a single cycle. Thus (ii) is true with \(\partial_G = \partial_L(L)\) and since \(\partial_V(\partial_0)\) is contained in the interior of \(\partial_L(L)\) and (b) holds, (iii) is also true. Here and henceforth we always obtain the outermost boundary using the same labelling procedure as described in the first paragraph.

In the rest we prove (a)-(b). By definition, if a cycle satisfies (a) and (b), it must be unique. Suppose not and there are two cycles \(A_1 = (x_1,...,x_t)\) and \(A_2 = (y_1,...,y_l)\) in \(G_B\) that satisfy (a) and (b). If \(A_1 \neq A_2,\) then there exists a vertex \(x_j\) of \(A_1\) not in \(A_2.\) Without loss of generality, suppose it belongs to the exterior of \(A_2.\) The corresponding square \(S_j\) with centre \(x_j\) is in \(\Lambda_0\) but lies in the exterior of \(A_2,\) contradicting the fact that \(A_2\) satisfies~(b).

To see (a) is true, let \(\Lambda_e \subset \Lambda_0\) be the set of vacant squares that share a vertex with some occupied square in \(C(0)\) and lie in the exterior of \(\partial_0.\) (By exterior of \(\partial_0,\) we mean exterior to every cycle in \(\partial_0\) and a square is in the interior of~\(\partial_0\) if it is in the interior of some cycle in~\(\partial_0.\)) If \(S\) is a square in \(C_V(\partial_0),\) then every corner of \(S\) either is the centre of a vacant square of~\(\Lambda_e\) or is the centre of an occupied square in \(C(0).\) This is true because every square that lies in the exterior of \(\partial_0\) and shares a vertex with~\(\partial_0\) is necessarily vacant. Also, \(\partial_V(\partial_0)\) contains \(\partial_0\) and therefore all occupied squares of \(C(0)\) in its interior. In particular, no point in the interior of \(\partial_0\) is in \(\partial_V(\partial_0).\) Thus, every vertex of \(\partial_V(\partial_0)\) is the centre of a vacant square in \(\Lambda_e.\)

%If some vertex in \(\partial_V(\partial_0)\) belonged to the centre of an occupied square, then

%and the edges of \(\partial_d\) are straight lines joining the centres of vacant squares.

%To see why (b) is true...

%We let \(\Lambda'_0 \subset \Lambda_0\) be the set of vacant squares that share a vertex with some occupied square in \(C(0)\) and lie in the exterior of \(\partial_0.\) By exterior of \(\partial_0,\) we mean exterior to every cycle in \(\partial_0\) and we say a square is in the interior of \(\partial_0\) if it is in the interior of some cycle in \(\partial_0.\) We claim that the centre of every square in \(\Lambda'_0\) is present in \(\partial_d.\)

To see (b) is true, we again use the fact that \(\partial_0\) is contained in the interior of~\(\partial_V(\partial_0)\) and thus every square in \(C(0)\) is contained in the interior of \(\partial_V(\partial_0).\) Also, any vacant square in \(\Lambda_0\setminus \Lambda_e\) that is contained in the interior of some cycle in~\(\partial_0\) is also contained in the interior of \(\partial_V(\partial_0).\) It only remains to see that no vacant square in \(\Lambda_e\) lies in the exterior of \(\partial_V(\partial_0).\) Suppose that a square in \(\Lambda_e\) lies in the exterior of \(\partial_V(\partial_0),\) then it is star adjacent to some occupied square \(S\) in \(C(0).\) But this means that some point in \(\partial_0\) either lies in the exterior of \(\partial_V(\partial_0)\) or intersects \(\partial_V(\partial_0),\) a contradiction since \(\partial_V(\partial_0)\) satisfies:
%Moreover, we have
\begin{equation} \label{eq1}
\text{Every point in } \partial_V(\partial_0) \text{ is at a distance of at least }0.5 \text{ from } \partial_0. \;\;\;\;\;\;\;\;
\end{equation}
This completes the proof of the theorem.

To see (\ref{eq1}) is true, let \(v_1\) and \(v_2\) be centres of two plus adjacent vacant squares \(W_1\) and \(W_2\) in \(\Lambda_e.\) Both \(v_1\) and \(v_2\) are at least a distance of \(0.5\) from any point in \(\partial_0.\) If the edge \(e\) common to \(W_1\) and \(W_2\) belongs to some cycle \(C\) in \(\partial_0,\) then either \(W_1\) or \(W_2\) is contained in the interior of \(C,\) a contradiction. Thus every point in the edge joining \(v_1\) and \(v_2\) is at a distance of at least \(0.5\) from \(\partial_0.\) Since \(\partial_V(\partial_0)\) consists only of edges joining centres of vacant squares in \(\Lambda_e,\) this proves (\ref{eq1}). \(\qed\)

\section{Proof of Theorem~\ref{thm5}}\label{pf5}
\emph{Proof of Theorem~\ref{thm5}}:
Let \(G^+_C\) denote the graph with vertex set as corners of squares of \(C^+(0)\) and edge set as edges of such squares. Let \(\partial^+_0 = (e_1,...,e_t)\) be the outermost boundary cycle in \(G^+_C\) for the component \(C^+(0)\) obtained from Theorem~2 of Ganesan (2015) and let \(v_i\) be the end vertex common to \(e_i\) and \(e_{i+1}\) for \(1 \leq i \leq t\) and \(v_t\) be the end vertex common to \(e_t\) and \(e_1.\)

%(show uniqueness later...after obtaining the outermost boundary...)

%(?? write later..

We first obtain the outermost boundary \(\partial_H\) in (ii) of the statement of the theorem, containing all vacant squares of \(\Lambda^+_0\) in its interior and then use \(\partial_H\) to obtain the corresponding star connected \(S-\)cycle \(H_{out}.\) Let \(\Lambda^+_1\) be the set of all vacant squares in \(\Lambda^+_0\) that lie in the exterior of \(\partial^+_0\) and let \(\Lambda^+_2 = \{W_i\}_{1 \leq i \leq t}\) be the squares in \(\Lambda^+_1\) such that \(W_i\) shares edge \(e_i\) with \(\partial^+_0.\) Let \(D_i\) denote the cycle formed by joining all squares in \(\Lambda^+_1\) that share an edge with \(W_i.\) We use the following properties below.
\begin{eqnarray}
&&\text{Every square in \(\Lambda^+_1\) is contained in some \(D_i, 1 \leq i \leq n.\) }\;\;\;\;\;\;\;\;\label{eqa}\\
&&\text{The cycle \(D_i\) contains the edge \(e_i.\)}\;\;\;\;\;\;\;\;\label{eqb}
\end{eqnarray}
%\begin{itemize}
%\item{\text{Every square in \(\Lambda^+_1\) shares an edge with some square in \(\Lambda^+_2.\)}}
%\item{\text{Every square in \(\Lambda^+_1\) shares an edge with some square in \(\Lambda^+_2.\)}}
%\end{itemize}
To see (\ref{eqa}) is true, we  note that every square in \(\Lambda^+_1\) shares an edge with some square in \(\Lambda^+_2.\) Suppose that \(S \in \Lambda^+_1\) shares only a vertex \(v_j\) with \(\partial_0^+.\) The square \(W_j\) is attached to edge \(e_j\) and necessarily shares an edge with \(S.\) To see (\ref{eqb}) is true, we note that \(e_i\) is adjacent to an occupied square and \(W_i\) and therefore cannot be in the interior of \(D_i,\) which only consists of vacant squares.

%(see fig...) (to write more..?? see crflly..)

To obtain \(\partial_H,\) we proceed iteratively and let \(C_0 = \partial_0^+ = (e_1,...,e_t).\) Before the iteration begins, all the vacant squares in \(D_1\) are contained in the exterior of \(\partial_0^+.\) Also the cycles \(D_1\) and \(C_0\) share the edge \(e_1\) and thus at least two points in common. We merge \(D_1\) and \(C_0\) using Theorem~3 of Ganesan (2015) to obtain a new cycle \(C_1\) that contains \(e_1\) in its interior. This is true because both the squares (one occupied and one vacant) that share \(e_1\) are now in the interior of \(C_1.\) At the first step of the iteration, we check if \(e_{2} \in C_{1}.\) If so, then since \(e_2\) is a boundary edge, it is attached an occupied square in the interior of \(C_{1}\) and the vacant square \(W_2\) in the exterior. We merge \(C_{1}\) and \(D_2\) and obtain a new cycle \(C_2.\) The edge \(e_2\) now belongs to the interior of \(C_2.\)

If \(e_2 \notin C_{1},\) it belongs to the interior of \(C_{1}.\) We then check if all squares of \(D_2\) are in the interior of \(C_{1}\) and if so, we set \(C_{1} = C_2\) and finish this iteration step. Else, let \(Y_1, Y_2\) and \(Y_3\) be the squares that are plus adjacent to \(W_2\) and share edges \(g_1, g_2\) and \(g_3,\) respectively with \(W_2\) and let \(X_2\) be the occupied square of \(C^+(0)\) that shares \(e_2\) with \(W_2.\) If \(Y_1\) lies in the exterior of~\(C_{1},\) then necessarily  \(g_1 \in C_{1}\) and thus has two vertices in common with \(C_{1}.\) If further \(Y_1 \in \Lambda^+_1,\) we merge \(Y_1\) with~\(C_{1}\) using Theorem~3 of Ganesan (2015) to obtain a new cycle \(C'_1.\) Using \(C'_1,\) we do the same for \(Y_2\) to get a new cycle \(C''_1\) and using \(C''_1\) and we do the above for \(Y_3\) to obtain a final cycle \(C_2.\) This completes the the iteration step  and in the next iteration step, we then repeat the above procedure with \(e_3\) and \(C_2.\)

%(see figure...)?? for Y_1, Y_2...to see crflly...

Continuing this way we get final a single cycle \(D_{fin}\) and by construction, it satisfies the following properties:\\
(a1) The cycle \(D_{fin}\) contains only edges of vacant squares in \(\Lambda^+_1.\)\\
(a2) Every edge of the outermost boundary \(\partial_0^+\) lies in the interior of \(D_{fin}.\)\\
(a3) For each \(i, 1 \leq i \leq t,\) the interior of \(D_i\) is contained in the interior of~\(D_{fin}.\)\\
From (a2) we know that \(\partial_0^+\) and therefore all squares contained in the interior of \(\partial_0^+,\) are also in the interior of \(D_{fin}.\) In addition, using (a3) and (\ref{eqa}), we have
\begin{equation}
\text{Every square in \(C^+(0) \cup \Lambda^+_0\) is contained in the interior of \(D_{fin}.\)}\label{eq_dfin2}
\end{equation}
The cycle \(D_{fin}\) is unique in the sense that if any cycle \(C\) satisfies (a1) and (\ref{eq_dfin2}), then \(C = D_{fin}.\)  To see this is true, suppose there is a
cycle \(C\) distinct from \(D_{fin} \) that satisfies (a1) and (\ref{eq_dfin2}) and suppose \(C\) contains an edge \(e\) in the exterior of \(D_{fin}.\) The edge \(e\) is adjacent to two squares \(Z_{e,1}\) and \(Z_{e,2}\) both of which lie in the exterior of \(D_{fin}.\) Moreover, at least one of \(Z_{e,1}\) or \(Z_{e,2}\) is a vacant square in \(\Lambda^+_1.\)  But \(D_{fin}\) satisfies (\ref{eq_dfin2}) and therefore \(e\) cannot be an edge of a square in
\(C^+(0) \cup \Lambda^+_0,\) a contradiction.

%To see that (\ref{eq_dfin}) is true, suppose after the first iteration, we have that \(e_2\) belongs to the interior of \(C_1.\)

%The cycle \(D_{fin}\) is the required \(\partial_d^+\) and satisfies (ii) and (iii). (and to see crflly...)
%Let \(\{Z_j\}_{1 \leq j \leq r}\) be the squares in \(\Lambda_0^+\) that share a vertex with \(D_{fin}.\) we obtain the cycle this way??

%The uniqueness above guarantees that we obtain the same cycle whatever be the order in which the edges of \(\partial^+_0\) are considered in the iteration above.

To see that the cycle \(D_{fin}\) is the required \(\partial_H,\) we now use \(D_{fin}\) to obtain the star connected \(S-\)cycle \(H_{out}.\) Let the cycle \(D_{fin} = (f_1,...,f_r).\) There exists a unique square \(Z_1 \in \Lambda^+_1\) that has edge \(f_1.\) Because, if two squares \(Z_1\) and \(Z'_1\) in \(\Lambda_1^+\) are attached to \(f_1,\) i.e. share the edge \(f_1,\) then necessarily one of them is in the exterior of \(D_{fin}.\) This contradicts (\ref{eq_dfin2}).

Similarly there exists a unique vacant square \(Z_2 \in \Lambda^+_1\) that has edge \(f_2.\) If \(Z_2 = Z_1,\) we proceed to \(f_3,\) else \(Z_2\) is star adjacent to \(Z_1\) and we add \(Z_2\) to the existing sequence and obtain \((Z_1,Z_2).\) Continuing this way, we obtain a final sequence of squares \(L_1 = (Z_1,...,Z_s)\) such that \(Z_i\) is star adjacent to \(Z_{i+1}\) for \(1 \leq i \leq s-1\) and \(Z_s\) is adjacent to \(Z_1.\) It only remains to see that this sequence \(L_1\) is the desired \(H_{out}.\)

By construction, the sequence \(L_1\) obtained is unique and (i) is also true. The following Lemma and (\ref{eq_dfin2}) imply that the sequence \(L_1\) also satisfies (ii) and (iii).
\begin{Lemma}\label{bd_l1} The outermost boundary \(\partial_0(L_1)\) of \(L_1\) is \(D_{fin}.\)
\end{Lemma}
This completes the proof of Theorem~\ref{thm5}.\(\qed\)

\emph{Proof of Lemma~\ref{bd_l1}}: Let \(G_L\) denote the connected graph containing vertices as corners of squares in \(\Lambda^+_1\) and edge set being  the edges of squares in \(\Lambda^+_1.\) From Theorem~1 of Ganesan (2015), we know that the outermost boundary \(\partial_0(L_1)\) is a connected union of cycles in \(G_L\) and every square in \(L_1\) is contained in the interior of some cycle in \(\partial_0(L_1).\) %As in the proof of Theorem~\ref{thm1} above, we label squares as \(0\) and \(1\) obtain the corresponding outermost boundary \(\partial_0(L_1).\) %by labelling all squares in \(L_1\) as \(1,\) all squares sharing a vertex with \(L_1\) as in Theorem~\ref{thm1} above.

By (\ref{eq_dfin2}), all squares of \(\Lambda^+_0\) are contained in the interior of \(D_{fin}.\) So every edge in \(\partial_0(L_1)\) either belongs to \(D_{fin}\) or is contained in its interior. If there exists an edge \(e\) of \(D_{fin}\) not in \(\partial_0(L_1),\) it necessarily lies in the exterior of all cycles in \(\partial_0(L_1).\) Also there exists some \(Z_j \in L_1\) that contains \(e\) as an edge. Thus \(Z_j\) lies in the exterior of \(\partial_0(L_1),\) a contradiction. \(\qed\)

\renewcommand{\theequation}{\thesection.\arabic{equation}}

\subsection*{Acknowledgement}
I thank Professor Rahul Roy for crucial comments and NISER for my fellowship.

\bibliographystyle{plain}

\end{document}